Definition and Existence of the Eigenderivative

Kerry M. Soileau

September 15, 2011

ABSTRACT

We define the eigenderivatives of a linear operator on any real or complex Banach space, and give a sufficient condition for their existence.

**MOTIVATION**

If $K \in \text{End}(X)$[1] is a linear operator on $X$, a real or complex Banach space, we hypothesize that the eigenvectors and eigenvalues of an operator "nearly equal" to $K$ will typically be "nearly equal" to those of $K$. In this paper we explore this question and give a definition of the "eigenderivative," i.e. the "rates of change" of eigenvectors and eigenvalues.

**DEFINITION OF EIGENDERIVATIVES**

Definition: Let $X$ be a real or complex Banach space, and let $K \in \text{End}(X)$ be a linear operator with a complete set of unit eigenvectors $\{e_i(K)\}_{i=1}^{N}$ with distinct eigenvalues $\{\lambda_i(K)\}_{i=1}^{N}$, where $N \leq \infty$. Let $F = \mathbb{C}$ or $\mathbb{R}$ identify the field over $X$. Suppose $J$ is a linear operator such that for

---

[1] The notation $\text{End}(X)$ denotes the linear operators mapping $X$ into itself.

each $i \in \{1,2,3,\cdots,N\}$ [2], $\sum_{\substack{j=1 \\ j \neq i}}^{N} \frac{j_{i,j}}{\lambda_i(K) - \lambda_j(K)} e_j(K)$ and $\sum_{\substack{j=1 \\ j \neq i}}^{N} \frac{j_{i,j}}{\lambda_i(K) - \lambda_j(K)} (J - j_{i,i}) e_j(K)$

converge, where $Je_i(K) = \sum_{j=1}^{N} j_{i,j} e_j(K)$. Then for each $i \in \{1,2,3,\cdots,N\}$ we define the

<u>eigenderivatives</u> $\Delta_i : \text{End}(X) \times \text{End}(X) \to X$ and $\Lambda_i : \text{End}(X) \times \text{End}(X) \to F$ according to

$$\Lambda_i(K,J) = j_{i,i} \text{ and } \Delta_i(K,J) = \sum_{\substack{j=1 \\ j \neq i}}^{N} \frac{j_{i,j}}{\lambda_i(K) - \lambda_j(K)} e_j(K).$$

These eigenderivatives verify

$$(K + hJ)(e_i(K) + h\Delta_i(K,J)) = (\lambda_i(K) + h\Lambda_i(K,J))(e_i(K) + h\Delta_i(K,J))$$
$$+ h^2 \sum_{\substack{j=1 \\ j \neq i}}^{N} \frac{j_{i,j}}{\lambda_i(K) - \lambda_j(K)} (J - j_{i,i}) e_j(K).$$

<u>Proposition 1</u>: If $J$ is bounded, and if for each $i \in \{1,2,3,\cdots,N\}$ we have that $\sup_j |j_{i,j}| < \infty$ and

$\sum_{\substack{j=1 \\ j \neq i}}^{N} \frac{1}{\lambda_i(K) - \lambda_j(K)}$ converges absolutely, then the eigenderivatives exist.

<u>Proof</u>:

$$\sum_{\substack{j=1 \\ j \neq i}}^{N} \left\| \frac{j_{i,n}}{\lambda_i(K) - \lambda_j(K)} e_j(K) \right\| \leq \sum_{\substack{j=1 \\ j \neq i}}^{N} \left| \frac{j_{i,n}}{\lambda_i(K) - \lambda_j(K)} \right| \|e_j(K)\| = \sum_{\substack{j=1 \\ j \neq i}}^{N} \left| \frac{j_{i,n}}{\lambda_i(K) - \lambda_j(K)} \right|$$
$$\leq \sum_{\substack{j=1 \\ j \neq i}}^{N} |j_{i,n}| \frac{1}{|\lambda_i(K) - \lambda_j(K)|} \leq \left( \sup_j |j_{i,j}| \right) \sum_{\substack{j=1 \\ j \neq i}}^{N} \frac{1}{|\lambda_i(K) - \lambda_j(K)|} < \infty$$

---

[2] In case $N = \infty$, $\{1,2,3,\cdots,N\}$ means the set of positive integers.

hence $\sum_{\substack{j=1 \\ j \neq i}}^{N} \dfrac{j_{i,n}}{\lambda_i(K) - \lambda_j(K)} e_j(K)$ converges.

$$\sum_{\substack{j=1 \\ j \neq i}}^{N} \left\| \dfrac{j_{i,j}}{\lambda_i(K) - \lambda_j(K)} (J - j_{i,i}) e_j(K) \right\| \leq \sum_{\substack{j=1 \\ j \neq i}}^{N} \left\| \dfrac{j_{i,j}}{\lambda_i(K) - \lambda_j(K)} (J e_j(K) - j_{i,i} e_j(K)) \right\|$$

$$\leq \sum_{\substack{j=1 \\ j \neq i}}^{N} \left\| \dfrac{j_{i,j}}{\lambda_i(K) - \lambda_j(K)} J e_j(K) \right\| + \sum_{\substack{j=1 \\ j \neq i}}^{N} \left\| \dfrac{j_{i,j} j_{i,i}}{\lambda_i(K) - \lambda_j(K)} e_j(K) \right\|$$

$$= \sum_{\substack{j=1 \\ j \neq i}}^{N} \left| \dfrac{j_{i,j}}{\lambda_i(K) - \lambda_j(K)} \right| \| J e_j(K) \| + \sum_{\substack{j=1 \\ j \neq i}}^{N} \left| \dfrac{j_{i,j} j_{i,i}}{\lambda_i(K) - \lambda_j(K)} \right| \| e_j(K) \|$$

$$= \sum_{\substack{j=1 \\ j \neq i}}^{N} |j_{i,j}| \dfrac{1}{|\lambda_i(K) - \lambda_j(K)|} \| J e_j(K) \| + |j_{i,i}| \sum_{\substack{j=1 \\ j \neq i}}^{N} |j_{i,j}| \dfrac{1}{|\lambda_i(K) - \lambda_j(K)|}$$

$$\leq \|J\| \sum_{\substack{j=1 \\ j \neq i}}^{N} |j_{i,j}| \dfrac{1}{|\lambda_i(K) - \lambda_j(K)|} + |j_{i,i}| \sum_{\substack{j=1 \\ j \neq i}}^{N} |j_{i,j}| \dfrac{1}{|\lambda_i(K) - \lambda_j(K)|}$$

$$\leq \|J\| \sup_j |j_{i,j}| \sum_{\substack{j=1 \\ j \neq i}}^{N} \dfrac{1}{|\lambda_i(K) - \lambda_j(K)|} + \left( \sup_j |j_{i,j}| \right)^2 \sum_{\substack{j=1 \\ j \neq i}}^{N} \dfrac{1}{|\lambda_i(K) - \lambda_j(K)|}$$

$$= \left( \|J\| + \sup_j |j_{i,j}| \right) \left( \sup_j |j_{i,j}| \right) \sum_{\substack{j=1 \\ j \neq i}}^{N} \dfrac{1}{|\lambda_i(K) - \lambda_j(K)|} < \infty$$

hence $\sum_{\substack{j=1 \\ j \neq i}}^{N} \dfrac{j_{i,j}}{\lambda_i(K) - \lambda_j(K)} (J - j_{i,i}) e_j(K)$ converges.

Examples: Let $X$ be the completion of the vector space spanned by the simple functions

$f_n : (0, \infty) \to \mathbb{R}$ verifying $f_n(x) = \begin{cases} 1 & n \leq x < n+1 \\ 0 & \text{elsewhere} \end{cases}$ for positive integers $n$, with $F = \mathbb{R}$. Let the

$\|\cdot\|$ be given by the $L^2$ norm, i.e. $\|f\| = \sqrt{\int_0^\infty f(x)^2 \, dx}$.

1. *An unbounded operator, for which eigenderivatives exist.*

Let $K$ satisfy $Kf_n(x) = nf_n(x)$. Then $\lambda_n(K) = n$ and $e_n(K) = f_n$. Let $J$ be the operator satisfying $Jf_n(x) = \sum_{k=1}^{\infty} \frac{1}{\sqrt{n+k}} f_k(x)$. Note that $\|Jf_n(x)\| = \infty$ thus $J$ is unbounded. Next note that

$$\sum_{\substack{j=1 \\ j\neq i}}^{\infty} \left\| \frac{j_{i,j}}{\lambda_i(K) - \lambda_j(K)} e_j(K) \right\| = \sum_{\substack{j=1 \\ j\neq i}}^{\infty} \left\| \frac{1}{i-j} \frac{1}{\sqrt{i+j}} f_j(x) \right\| = \sum_{\substack{j=1 \\ j\neq i}}^{\infty} \frac{1}{|i-j|} \frac{1}{\sqrt{i+j}} < \infty$$

Next,

$$\left\| \sum_{\substack{j=1 \\ j\neq i}}^{\infty} \frac{j_{i,j}}{\lambda_i(K) - \lambda_j(K)} (J - j_{i,i}) e_j(K) \right\|$$

$$= \left\| \sum_{\substack{j=1 \\ j\neq i}}^{\infty} \left( \frac{1}{\sqrt{2i}} \frac{1}{j-i} \frac{1}{\sqrt{i+j}} - \sum_{\substack{k=1 \\ k\neq i}}^{\infty} \frac{1}{k+j} \frac{1}{k-i} \frac{1}{\sqrt{i+k}} \right) f_j(x) - \sum_{\substack{k=1 \\ k\neq i}}^{\infty} \frac{1}{k^2 - i^2} \frac{1}{\sqrt{i+k}} f_i(x) \right\|$$

$$= \sqrt{\sum_{\substack{j=1 \\ j\neq i}}^{\infty} \left( \frac{1}{\sqrt{2i}} \frac{1}{j-i} \frac{1}{\sqrt{i+j}} - \sum_{\substack{k=1 \\ k\neq i}}^{\infty} \frac{1}{k+j} \frac{1}{k-i} \frac{1}{\sqrt{i+k}} \right)^2 + \sum_{\substack{k=1 \\ k\neq i}}^{\infty} \frac{1}{(k^2-i^2)^2} \frac{1}{i+k}} < \infty$$

Thus the requirements of the Definition are satisfied, and the eigenderivatives exist.

They are $\Lambda_i(K,J) = \frac{1}{\sqrt{2i}}$ and $\Delta_i(K,J) = \sum_{\substack{j=1 \\ j\neq i}}^{\infty} \frac{1}{(i-j)\sqrt{i+j}} f_j(x)$. In this case we have $\|\Delta_i(K,J)\| = \sqrt{\sum_{\substack{j=1 \\ j\neq i}}^{\infty} \frac{1}{(j-i)^2 (i+j)}}$, depicted in the following:

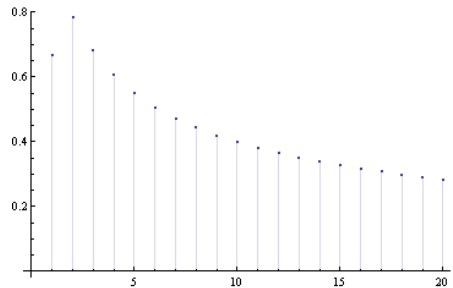

2. *A bounded operator, for which eigenderivatives exist.*

Let $K$ satisfy $Kf_n(x) = \frac{1}{n} f_n(x)$. Then $\lambda_n(K) = \frac{1}{n}$ and $e_n(K) = f_n$. Let $J$ be the operator satisfying $Jf_n(x) = \sum_{k=1}^{\infty} \frac{1}{n+k} f_k(x)$. Note that

$$\sup_{1 \leq n} \|Jf_n(x)\| = \sup_{1 \leq n} \sqrt{\sum_{k=1}^{\infty} \frac{1}{(n+k)^2}} < \frac{\pi^2}{6}$$

thus $J$ is bounded. Next note that

$$\left\| \sum_{\substack{j=1 \\ j \neq i}}^{\infty} \frac{j_{i,j}}{\lambda_i(K) - \lambda_j(K)} e_j(K) \right\| = i \sqrt{\sum_{\substack{j=1 \\ j \neq i}}^{\infty} \frac{j^2}{(j^2 - i^2)^2}} < \infty \text{ and}$$

$$\left\| \sum_{\substack{j=1 \\ j \neq i}}^{\infty} \frac{j_{i,j}}{\lambda_i(K) - \lambda_j(K)} (J - j_{i,i}) e_j(K) \right\| = \sum_{\substack{j=1 \\ j \neq i}}^{\infty} \left( \frac{1}{4} \frac{1}{(j+i)^2} + \left( \frac{ij}{j^2 - i^2} \right)^2 \sum_{\substack{k=1 \\ k \neq j}}^{\infty} \frac{1}{(j+k)^2} \right) < \infty.$$

Thus the requirements of the Definition are satisfied, and the eigenderivatives exist.

They are $\Lambda_i(K, J) = \frac{1}{2i}$ and $\Delta_i(K, J) = i \sum_{\substack{j=1 \\ j \neq i}}^{\infty} \frac{j}{j^2 - i^2} f_j(x)$. In this case we have

$$\|\Delta_i(K, J)\| = i \sqrt{\sum_{\substack{j=1 \\ j \neq i}}^{\infty} \frac{j^2}{(j^2 - i^2)^2}} \approx 0.907i.$$

**AN EXISTENCE CRITERION ON A HILBERT SPACE**

Proposition 2: Let $X$ be a real or complex Hilbert space, and let $K \in \text{End}(X)$ be a linear operator with a complete set of eigenvectors $\{e_i(K)\}_{i=1}^{N}$ with distinct eigenvalues $\{\lambda_i(K)\}_{i=1}^{N}$, where $N \leq \infty$. For any bounded linear operator $J \in \text{End}(X)$ and eigenvector $e_i(K)$, if for every $i \in \{1,2,3,\cdots,N\}$, we have that $\sum_{\substack{j=1 \\ j \neq i}}^{N} \frac{1}{\lambda_i(K) - \lambda_j(K)}$ converges absolutely, then the eigenderivatives exist and are given by $\Lambda_i(K,J) = \langle Je_i(K), e_i(K) \rangle$ and

$$\Delta_i(K,J) = \sum_{\substack{j=1 \\ j \neq i}}^{N} \frac{\langle Je_i(K), e_j(K) \rangle}{\lambda_i(K) - \lambda_j(K)} e_j(K).$$

Proof: For ease of notation this proof will employ the abbreviations $e_n = e_n(K)$ and $\lambda_n = \lambda_n(K)$. Fix $i \in \{1,2,3,\cdots,N\}$. Then $Je_i = \sum_{j=1}^{N} j_{i,j} e_j$ where $j_{i,j} = \langle Je_i, e_j \rangle$. This implies

$\sup_j |j_{i,j}| = \sup_j |\langle Je_i, e_j \rangle| \leq \sup_j \|J\| \|e_i\| \|e_j\| = \sup_j \|J\| = \|J\| < \infty$. Next, since by assumption $\sum_{\substack{j=1 \\ j \neq i}}^{N} \frac{1}{\lambda_i - \lambda_j}$ converges absolutely, by Proposition 1 we have that then the eigenderivatives exist, and are given by $\Lambda_i(K,J) = \langle Je_i, e_i \rangle$ and $\Delta_i(K,J) = \sum_{\substack{j=1 \\ j \neq i}}^{N} \frac{\langle Je_i, e_j \rangle}{\lambda_i - \lambda_j} e_j$

Positronic Software
League City, TX 77573

E-mail address: ksoileau@yahoo.com